\documentclass[english]{article}
\usepackage[T1]{fontenc}
\usepackage[latin9]{inputenc}
\usepackage{amsthm}
\usepackage{amstext}
\usepackage{amssymb}

\makeatletter

\newcommand{\lyxmathsym}[1]{\ifmmode\begingroup\def\b@ld{bold}
  \text{\ifx\math@version\b@ld\bfseries\fi#1}\endgroup\else#1\fi}

\newcommand{\lyxaddress}[1]{
\par {\raggedright #1
\vspace{1.4em}
\noindent\par}
}
\theoremstyle{plain}
\newtheorem{thm}{\protect\theoremname}
\ifx\proof\undefined
\newenvironment{proof}[1][\protect\proofname]{\par
\normalfont\topsep6\p@\@plus6\p@\relax
\trivlist
\itemindent\parindent
\item[\hskip\labelsep
\scshape
#1]\ignorespaces
}{%
\endtrivlist\@endpefalse
}
\providecommand{\proofname}{Proof}
\fi

\makeatother

\usepackage{babel}
\providecommand{\theoremname}{Theorem}

\begin{document}

\title{Product of Two Consecutive Fibonacci or Lucas Numbers Divisible by
their Prime Sum of Indices}

\author{Vladimir Pletser}

\maketitle

\lyxaddress{European Space Research and Technology Centre, ESA-ESTEC P.O. Box
299, NL-2200 AG Noordwijk, The Netherlands E-mail: Vladimir.Pletser@esa.int }
\begin{abstract}
{\normalsize{We show that the product of two consecutive Fibonacci
(respectively Lucas) numbers is divisible by the sum of their indices
if this sum is a prime number different from $5$ and in the form
$(4r+1)$ (respectively $(4r+3)$).}}{\normalsize \par}

\end{abstract}

\section{Introduction}

One of the most interesting divisibility properties of the Fibonacci
numbers is that for each prime $p$, there is a Fibonacci number $F_{n}$
such that $p$ divides $F_{n}$ (see, e.g. \cite{key-2}). More specifically,
for $p\neq5$, $p$ divides either $F_{p-1}$ if $p\equiv\pm1(mod\,5)$,
or $F_{p+1}$ if $p\equiv\pm2(mod\,5)$. For $p=5$, one has of course
$p=F_{p}$.

\section{Theorem}

Although already demonstrated differently in \cite{key-5,key-6},
a new demonstration of the following theorem is proposed in this paper.
\begin{thm}
If $p$ is prime and $r\in\mathbb{Z}^{+}$, 
\end{thm}
\begin{eqnarray}
p & = & \left(4r+1\right)\textnormal{ \textit{divides the product} }F_{2r}F_{2r+1},\textnormal{ \textit{except for} }p=5\label{eq:1}\\
p & = & \left(4r+3\right)\textnormal{ \textit{divides the product} }L_{2r+1}L_{2r+2}\label{eq:2}
\end{eqnarray}

\begin{proof}
For $p$ prime and $r,s,n,m\in\mathbb{Z}^{+}$, for odd primes $p=2s+1$,
one has 
\begin{equation}
L_{2s+1}-1=L_{2s+1}-L_{1}\label{eq:3}
\end{equation}
The transformations 
\begin{eqnarray}
L_{n+m}-\left(-1\right)^{m}L_{n-m} & = & 5F_{m}F_{n}\label{eq:4}\\
L_{n+m}+\left(-1\right)^{m}L_{n-m} & = & L_{m}L_{n}\label{eq:5}
\end{eqnarray}
(relations (17 a, b) in \cite{key-3} and relations (11) and (23)
in \cite{key-4}) can be used. 

\noindent (i) First, let $s$ be even, $s=2r$. Relation (\ref{eq:3})
yields respectively from (\ref{eq:4}) and (\ref{eq:5}), with $m=2r$
and $n=2r+1$,
\begin{eqnarray}
L_{4r+1}-1 & = & 5F_{2r}F_{2r+1}\label{eq:6}\\
L_{4r+1}+1 & = & L_{2r}L_{2r+1}\label{eq:7}
\end{eqnarray}
If $p=4r+1\neq5$ is prime, then either $p$ divides $F_{4r}$ if
$p\equiv\pm1(mod\,5)=29,41,61,\lyxmathsym{\ldots}$, or $p$ divides
$F_{4r+2}$ if $p\equiv\pm2(mod\,5)=13,17,37,\lyxmathsym{\ldots}$ 

\noindent On the other hand, one has (relation (13) in \cite{key-3})
\begin{eqnarray}
F_{4r} & = & F_{2r}L_{2r}\label{eq:8}\\
F_{4r+2} & = & F_{2r+1}L_{2r+1}\label{eq:9}
\end{eqnarray}
Let first $p\equiv\pm1(mod\,5)$, then $p$ divides $F_{4r}$ and
therefore from (\ref{eq:8}) also either $F_{2r}$ or $L_{2r}$. But
$p$ cannot divide $L_{2r}$. Let us assume the contrary. Suppose
that $p$ divides $L_{2r}$ and also $\left(L_{4r+1}-1\right)$, as
$L_{p}\equiv1\left(mod\, p\right)$ (see e.g. \cite{key-1}, \cite{key-7}).
It would mean from (\ref{eq:7}) that $p$ should also divide simultaneously
$\left(L_{4r+1}+1\right)$ which makes no sense. Therefore $p$ divides
$F_{2r}$ and not $L_{2r}$, and also $\left(L_{4r+1}\text{\textendash}1\right)$. 

\noindent The other case where $p\equiv\pm2(mod\,5)$ divides $F_{4r+2}$
is treated similarly. 

\noindent This means that all primes $p=4r+1\neq5$ divide the product
of two consecutive Fibonacci numbers of indices $2r$ and $2r+1$.
More precisely, if $p\equiv\pm1(mod\,5)=29,41,61,\lyxmathsym{\ldots}$,
then $p$ divides $F_{2r}$; if $p\equiv\pm2(mod\,5)=13,17,37,\lyxmathsym{\ldots}$,
then $p$ divides $F_{2r+1}$. 

\noindent (ii) Second, let $s$ be odd, $s=2r+1$. Relation (\ref{eq:3})
yields respectively from (\ref{eq:4}) and (\ref{eq:5}), with $m=2r+1$
and $n=2r+2$,
\begin{eqnarray}
L_{4r+3}+1 & = & 5F_{2r+1}F_{2r+2}\label{eq:10}\\
L_{4r+3}-1 & = & L_{2r+1}L_{2r+2}\label{eq:11}
\end{eqnarray}
If $p=4r+3$ is prime, then $p$ divides $F_{4r+2}$ if $p\equiv\pm1(mod\,5)=11,19,31,\lyxmathsym{\ldots}$;
or $p$ divides $F_{4r+4}$ if $p\equiv\pm2(mod\,5)=3,7,23,43,\lyxmathsym{\ldots}$
One has also 

\noindent 
\begin{eqnarray}
F_{4r+2} & = & F_{2r+1}L_{2r+1}\label{eq:12}\\
F_{4r+4} & = & F_{2r+2}L_{2r+2}\label{eq:13}
\end{eqnarray}
Like above, let first $p\equiv\pm1(mod\,5)$. Then $p$ divides $F_{4r+2}$
and therefore, from (\ref{eq:12}), also either $F_{2r+1}$ or $L_{2r+1}$.
But $p$ cannot divide $F_{2r+1}$. Let us assume the contrary. Suppose
that $p$ divides $F_{2r+1}$ and also $\left(L_{4r+3}-1\right)$.
It would mean from (\ref{eq:10}) that $p$ should also divide simultaneously
$\left(L_{4r+3}\text{+}1\right)$ which makes no sense. Therefore
$p$ divides $L_{2r+1}$ and not $F_{2r+1}$, and also $\left(L_{4r+3}\text{\textendash}1\right)$.
The other case where $p\equiv\pm2(mod\,5)$ divides $F_{4r+4}$ is
also treated similarly. 

\noindent This means that all primes $p=4r+3$ divide the product
of two consecutive Lucas numbers of indices $2r+1$ and $2r+2$. More
precisely, if $p\equiv\pm1(mod5)=11,19,31,\lyxmathsym{\ldots}$, then
$p$ divides $L_{2r+1}$; if$p\equiv\pm2(mod\,5)=3,7,23,43,\lyxmathsym{\ldots}$,
then $p$ divides $L_{2r+2}$. 

\noindent On the other hand, for $p=5$, one has obviously $L_{5}=11\equiv1(mod\,5)$. 
\end{proof}

\section{Acknowledgment }

Dr C. Thiel is acknowledged for helpful discussions.


\begin{thebibliography}{1}
\bibitem[1]{key-1} P. S. Bruckman, Lucas Pseudoprimes are Odd, Fibonacci
Quarterly 32, 155-157, 1994.

\bibitem[2]{key-4} R. A. Dunlap, The Golden Ratio and Fibonacci Numbers,
World Scientific Press, 1997. 

\bibitem[3]{key-7} T. Koshy, Fibonacci and Lucas Numbers with Applications,
John Wiley, New York, p. 410, 2001.

\bibitem[4]{key-6} J. Seibert, Fibonacci and Lucas Products Modulo
A Prime, Solution Problem B-1037, Fibonacci Quarterly, Vol. 46-47,
p. 88, 2008-2009

\bibitem[5]{key-2} M. R. Schroeder, Number Theory in Science and
Communication, 2nd edition, Springer-Verlag, 1986, 72-73. 

\bibitem[6]{key-3} S. Vajda, Fibonacci and Lucas Numbers, and The
Golden Section: Theory and Applications, Halsted Press, 1989. 

\bibitem[7]{key-5} H. C. Williams, Edouard Lucas and primality testing,
Canadian Math. Soc. Monographs 22, Wiley, New York, 1998.\end{thebibliography}
\end{document}